\documentclass[a4paper,12pt]{article}
\usepackage{amsmath}
\usepackage{amssymb}
\usepackage{amsthm}
\usepackage[utf8]{inputenc}

\newtheorem{theorem}{Theorem}

\title{Improved Lower Bounds  for Property B}
\author{Karl Grill \and Daniel Linzmayer%
\\Institute of Statistics and Mathematical Methods in Economics\\TU Wien}
\begin{document}
\maketitle
\begin{abstract}
If an $n$-uniform hypergraph can be 2-colored, then it is 
said to have property B. Erdős (1963) was the
first to give lower and upper bounds for
 the minimal size $\mathbf m(n)$ of an $n$-uniform  hypergraph without
property B. His asymptotic upper bound $O(n^22^n)$ still is the best
we know, his lower bound $2^{n-1}$ has seen a number of 
improvements, with the current best 
$\Omega$ $(2^n\sqrt{n/\log(n)})$ 
established by  Radhakrishnan 
and Srinivasan (2000).
Cherkashin and Kozik (2015) provided a simplified proof
of this result, using Pluhár's (2009) idea of a random 
greedy coloring. In the present paper, we use a refined version of
this argument to obtain improved lower bounds on $\mathbf m(n)$ for small values of $n$.
We also study $\mathbf m(n,v)$, the size of the smallest $n$-hypergraph without
property B having $v$ vertices.
\end{abstract}

\section{Introduction}

We consider an $n$-uniform hypergraph $H=(V,E)$ with vertex set $V$ of
cardinality $|V|=v$ and edge set 
\[E\subseteq \{A\subseteq V:|A|=n\}.\]
$H$ is said to have property B if it is  2-colorable, i.e., if its vertices
can be colored with two colors (traditionally called ``red'' and ``blue'')
in such a way that no edge is monochromatic. 

We let $\mathbf  m(n)$ denote the smallest number of edges that an $n$-uniform 
hypergraph without property B must have, and $\mathbf m(n,v)$ the smallest number
of edges in a hypergraph with $v$ vertices that does not have property B.

The exact values of $\mathbf m(n)$ are only known for $n = 1, 2, 3, 4$ so far. 
For $n=5$ and higher only bounds for the true solution are 
known\cite{oester}. Brute force calculation quickly reaches its limits due to 
the quickly increasing complexity of the problem for increasing $n$.

Erdős and Hajnal \cite{erdhaj} presented the first upper bound
$m(p)\le \genfrac(){0pt}1{2n-1}{n}.$
Erdős \cite{erdos}, using the probabilistic method, obtained the bounds
\begin{equation}\label{erdosupper}
2^{n-1}\le \mathbf m(n) \leq (1 + o(1)) e \log(2) n^2 2^{n-2}.\end{equation}
His upper bound is still the best asymptotic result available, though some
improvements have been made for small $n$ by constructive means.

In \cite{combprob3}, for even $v=O(n)$, he proves the bounds

\begin{equation}\label{erdosbounds}
\frac{\genfrac(){0pt}1{v}{n}}{2\genfrac(){0pt}1{v/2}{n}}\le \mathbf m(n,v)
\le 2v
\frac{\genfrac(){0pt}1{v}{n}}{2\genfrac(){0pt}1{v/2}{n}}.\end{equation}
This actually holds for any even $v>2n$, but for $v>n^2/2$ this is worse than
(\ref{erdosupper}), which is obtained for $v=n^2/2$.
So, for any $v$, we have upper and lower bounds for $\mathbf m(n,v)$ 
that differ only by a
factor of order $O(n^2)$. This is still quite big, but its polynomial
growth is slower than the exponential growth of the bounds on $\mathbf m(n,v)$,
a fact that will be important in our later considerations.

As for the lower bound, 
Goldberg and Russell\cite{gold} observed that for the smallest $v$ 
with $\mathbf m(n,v)\le m$, a hypergraph with $v$ vertices
achieving this bound must have the property
that every pair of vertices is contained in some edge, which, by 
Schoenheim's bound, implies
\begin{equation}\label{schoen}
m\ge\lceil\frac{v}{n}\lceil\frac{v-1}{n-1}\rceil\rceil.\end{equation}
Together with the random coloring bound
\begin{equation}\label{goldlow}\mathbf m(n,v)\ge\frac{\genfrac(){0pt}1{v}{n}}
{\genfrac(){0pt}1{\lfloor v/2\rfloor}{n}+
\genfrac(){0pt}1{\lceil v/2\rceil}{n}},\end{equation}
this gives
\begin{equation}\label{goldbound}\mathbf m(n)\ge\min_v\max\left(
\lceil\frac{v}{n}\lceil\frac{v-1}{n-1}\rceil\rceil,
\left\lceil\frac{\genfrac(){0pt}1{v}{n}}
{\genfrac(){0pt}1{\lfloor v/2\rfloor}{n}+
\genfrac(){0pt}1{\lceil v/2\rceil}{n}}
\right\rceil\right).\end{equation}
Beck \cite{beck1,beck2} used a recoloring of a random coloring and
succeeded in proving $2^{-n}\mathbf m(n)\to\infty$. Later
Spencer \cite{spencer} strengthened and simplified Beck's argument.
This idea was carried further
by
Radhakrishnan and Srinivasan \cite{srini}, yielding the best
lower bound currently known
\[\mathbf m(n)=\Omega(
2^n \sqrt{n/\log(n)}).\]
Cherkashin and Kozik \cite{cherk} gave a simpler proof of this
result, 
using the greedy coloring approach introduced by 
Pluhár \cite{pluhar}. We study this approach in more detail in the next 
section.

Many of these results generalize to more than two colors. The survey 
article by Raigorodskii and Cherkashin \cite{raicher} gives an account of 
various results along with their proofs.

\section{Greedy Coloring}

Pluhár \cite{pluhar} introduced a greedy coloring procedure: one starts with all
vertices red and arranges them in random order. Then one looks at the vertices
in sequence, changing the color of a vertex to blue if it is the last one in 
an otherwise all-red edge. 

By the nature of this algorithm, no monochromatic red edge can occur, and it is
obvious that for a two-colorable hypergraph  there is an ordering of the
vertices for which the greedy algorithm yields a proper coloring.
As the dependence of this procedure on the number of vertices is a bit of a
nuisance, Pluhár \cite{pluhar} introduced the notion of assigning independent, 
uniformly $[0,1]$
distributed weights to the vertices and arranging them in
increasing order of their weights. 

Using this idea, he obtained a simpler proof of the fact that $\mathbf m(n)2^{-n}\to
\infty$, although his result was weaker than that of Radhakrishnan and
Srinivasan.
The next step was performed by Cherkashin and Kozik \cite{cherk}. 
They utilized the random greedy coloring  method 
to construct a simpler proof for Radhakrishnan and Srinivasan's 
asymptotic 
result.

The central idea is the following: greedy coloring can only fail if it
produces a blue edge. By the nature of the coloring procedure, 
the first vertex in this edge must be the
last vertex in some other (otherwise red) edge. 
Thus, the probability that coloring fails can be estimated above by the
probability that some vertex is the last one in some edge and the first in some
other at the same time. We call such a vertex {\em critical}. Given that there
is a vertex with weight $x$, we can bound the conditional probability that
it is critical by each of the three probabilities that it is the last vertex
in some edge, that it is the first in some edge, or that both hold at the
same time. This yields the estimate
\begin{equation}\label{upperweight}\mathbb P(\mbox{Greedy coloring fails})\le 
\int_0^1\min(mnx^{n-1},mn(1-x)^{n-1},\gamma x^{n-1}(1-x)^{n-1})dx,
\end{equation}
where $m$ is the number of edges, and $\gamma$ is the number of 
ordered pairs of edges that have exactly one vertex in common. 
In some cases, one can find good estimates for $\gamma$, but most of
the time we have to make do with 
the trivial bound $\gamma\le m(m-1)$. 
It is easily seen that the right-hand side of (\ref{upperweight}) is
the minimum of
\begin{equation}\label{ckbound}
2mx^n+\gamma\int_x^{1-x}u^{n-1}(1-u)^{n-1}du\end{equation}
If this is less than $1$, then there is a positive probability that
greedy coloring succeeds, and so $\mathbf m(n)>m$.
Cherkashin and Kozik \cite{cherk} weaken this by applying the inequalities
$u(1-u)\le1/4$ and $\gamma\le m^2$ to obtain Radhakrishnan and Srinivasan's
\cite{srini} result; Aglave et al.\ \cite{shanni} reduce the case $n=5$,
$m=28$  to
$v=23$, as all other values of $v$ are ruled out by (\ref{schoen}) and
(\ref{goldlow}). 
They improve the upper
bound on $\gamma$ to $\gamma\le 670$, which is enough to give a value less
than $1$
in equation (\ref{ckbound}), proving $\mathbf m(5)\ge 29$. 

We go back to considering $\mathbf m(n,v)$ with its explicit dependence on the
number of vertices $v$. 
As a consequence, instead of the continuous distribution of the weights, 
we can work directly with the discrete distribution of the random permutation. 
By the same reasoning that led us to (\ref{upperweight}), we can get an 
upper bound for
the probability $p$ that there is a critical vertex:
\[p\le \sum_{k=0}^v p(k),\]
where $p(k)$ denotes the probability that the vertex in position $k$ is
critical; this can be estimated above by
\begin{equation}\label{upperperm}
\sum_{k=0}^v
p(k)\le\sum_{k=0}^v
 \frac{1}{v}\min\left(mn\frac{
\genfrac(){0pt}1{k-1}{n-1}
}{
\genfrac(){0pt}1{v-1}{n-1}
},
mn\frac{
\genfrac(){0pt}1{v-k}{n-1}
}{
\genfrac(){0pt}1{v-1}{n-1}
},
\gamma\frac{
\genfrac(){0pt}1{k-1}{n-1}
\genfrac(){0pt}1{v-k}{n-1}
}{
\genfrac(){0pt}1{v-1}{n-1}
\genfrac(){0pt}1{v-n}{n-1}
}\right).\end{equation}

It may be worth noting that for $v\to\infty$, the right-hand side of
(\ref{upperperm}) converges to (\ref{upperweight}) and, of course, the three
terms in the minimum are bounds for the probabilities that $k$ is the
last vertex in some edge, the first in some edge, or both.

In equation (\ref{upperperm}), the first term in the minimum is smaller
than the second for $k<\lfloor (v+1)/2\rfloor$, and if the smaller of these
is not greater than the third, then the right-hand side evaluates as

\[m\frac{\genfrac(){0pt}1{\lfloor v/2\rfloor}{n}
+\genfrac(){0pt}1{\lceil v/2\rceil}{n}}{\genfrac(){0pt}1{v}{n}},\]
so we can only get an improvement over (\ref{goldlow}), if, for some
$k<(v+1)/2$, we have
\[\gamma \genfrac (){0pt}0{v-k}{n-1}<mn\genfrac(){0pt}0{v-n}{n-1}.\]
As long as we do not have a better estimate than $\gamma\le m(m-1)$, 
(\ref{erdosbounds}) implies that we need 
\begin{equation}\label{vbound}
\genfrac(){0pt}0{v-1}{n-1}<n\genfrac(){0pt}0{v-n}{n-1},\end{equation}
which in turn implies $v>\frac{n^2}{\log(n)}(1+o(1))$. 

We can make some slight improvements to our estimate of $\gamma$: on one hand,
if we know that a certain pair of vertices is contained in $r$ edges, then
obviously
\[\gamma\le m(m-1)-r(r-1),\]
and another upper bound is obtained by counting the pairs that have a certain
vertex in common: let $l_i$ denote the number of occurrences of vertex $i$,
and $r_i$ the maximum frequency of a pair $\{i,j\}, j\neq i$. Then
\[\gamma\le\sum_{i=1}^v\left(l_i(l_i-1)-r_i(r_i-1)\right).\]

Using these ideas, we calculated lower bounds, letting $n$ vary in the range 
$5$ to $9$, and $v$ from $2n+1$ to $200$. In fact, in tune with our
considerations in (\ref{vbound}), we observe that, for small
values of $v$, our bound and the Goldberg-Russell bound (\ref{goldbound})
agree. For $n=5$, we get
improved lower bounds for $17\le v\le 25$.

For these calculations and those mentioned below, we use small C programs. As
all the probabilities in question, apart from the continuous Cherkashin-Kozik
bound (\ref{ckbound}) are rational, we can perform exact calculations using
an appropriate bignum library, we chose GNU {\tt libgmp}. Multiplying with
a common denominator, we can work with integers, which yields a considerable
improvement of efficiency over rational number computations.

For other values of $n$, too, this procedure
gives us improved lower
bounds for $\mathbf m(n)$. We summarize these bounds for small values of $n$
 together with those
obtained from Cherkashin and Kozik's \cite{cherk} continuous procedure,
that is independent of $v$, and the basic estimate by
Goldberg and Russell \cite{gold}
in table~\ref{table1}. For the sake of completeness, we have also included
the best known upper bounds as reported in \cite{shanni}.

\begin{table}[ht]
\begin{center}
\begin{tabular}{|l|c|c|c|c|c|}
\hline
$n$ & $5$ & $6$ & $7$ & $8$ & $9$ \\\hline
Goldberg and Russell \cite{gold} & 28 & 51 & 94 & 174 & 328 \\\hline
Cherkashin and Kozik \cite{cherk} & 27(29*) & 57 & 119 & 248 & 516\\\hline
Discrete greedy eq.\ (\ref{upperperm}) & 30 & 62 & 126 & 259 & 533\\\hline
Current upper bound \cite{shanni} & 51 & 147 & 421 & 1212 & 2401\\\hline
\end{tabular}
\caption{First lower bounds\\{\small *The lower bound 29 obtained in Aglave et 
al.\cite{shanni} improves the basic Cherkashin-Kozik estimate by giving a 
tighter estimate of the number of critical pairs and thusly goes beyond the 
basic procedure considered here; for the sake of completeness, however, we 
decided to add it in parentheses}}
\label{table1}
\end{center}
\end{table}

It should be mentioned that, although we get improved lower bounds for some
finite values of $n$, this approach
does not change the asymptotic result. As a matter 
of fact, the sum in (\ref{upperperm}), properly normalized, converges to the
integral in (\ref{upperweight}).

\section{Locking a vertex in place}
We consider a particular pair of edges that have only one vertex in common.
The probability that this pair becomes critical with their common point
in position $k$ evaluates as
\[\frac{\genfrac(){0pt}1{k-1}{n-1}\genfrac(){0pt}1{v-k}{n-1}}
{v\genfrac(){0pt}1{v-1}{n-1}\genfrac(){0pt}1{v-n}{n-1}}.\]

This takes its maximum for $k=\lceil v/2\rceil$. Thus, it seems plausible
that we might be able to improve our estimates by locking a vertex with a
small degree in place there. So, let $l$ denote the smallest degree
of a vertex, and let us put the associated vertex in position
$v_1=\lceil v/2\rceil$. There are $l$ edges that contain the selected
vertex, and $m-l$ that don't. By the pigeon-hole principle, we know
that
\[\frac{v-1}{n-1}\le l\le \frac{mn}{v}:\]
for the lower bound, remember that we assume that all pairs $\{i,j\}$ with
$j\neq i$ be contained in some edge, so $l(n-1)\ge v-1$, and on the other
hand, the sum of all $v$ vertex degrees is $mn$, so $lv\le mn$.

This assumption changes our bounds for the various probabilities, e.g.,
for $k<v_1$ we have
\[p_1(k)=n(m-l)\frac{\genfrac(){0pt}1{k-1}{n-1}}
{(v-1)\genfrac(){0pt}1{v-2}{n-1}},\]
\[p_2(k)=n(m-l)\frac{\genfrac(){0pt}1{v-k-1}{n-1}}
{(v-1)\genfrac(){0pt}1{v-2}{n-1}}+
nl\frac{\genfrac(){0pt}1{v-k-1}{n-2}}
{(v-1)\genfrac(){0pt}1{v-2}{n-2}},\]
\[p_3(k)=(m-l)(m-l-1)\frac{
\genfrac(){0pt}1{k-1}{n-1}\genfrac(){0pt}1{v-k-1}{n-1}}
{(v-1)\genfrac(){0pt}1{v-2}{n-1}\genfrac(){0pt}1{v-n-1}{n-1}}+
(m-l)l\frac{
\genfrac(){0pt}1{k-1}{n-1}\genfrac(){0pt}1{v-k-1}{n-2}}
{(v-1)\genfrac(){0pt}1{v-2}{n-1}\genfrac(){0pt}1{v-n-1}{n-2}}\]
as upper bounds for the probabilities that the vertex in position $k$ is
last in some edge, first in some edge, or both.

As it turns out, this approach serves to improve our lower bounds to
$\mathbf m(5)\ge 31$, $\mathbf m(6)\ge 63$, $\mathbf m(7)\ge 127$, $\mathbf m(8)\ge 261$ and
$\mathbf m(9)\ge 537$.

We can go further and fix more than one vertex in the center spots.
 Unfortunately, this quickly gets complicated, as we have
to control the numbers of edges that contain a particular subset of the
selected vertices, which amounts to $2^s-1$ numbers for $s$ selected 
vertices. In the case $s=2$, we need to control the
numbers $l_1$, $l_2$ and $l_{12}$ of edges containing only the first
selected vertex, only the second, or both. This 
 is still simple enough to carry out the
calculations for all combinations of $n$ and $v$ that we considered 
in the previous section. 
In this way, we obtain the improvements
$\mathbf m(7)\ge 128$ and  $\mathbf m(8)\ge 262$.

The case $s=3$ affords working with the seven numbers $l_A$ of edges
containing exactly the (non-empty) subset  $A$ of $\{1,2,3\}$. We apply this 
 to selected combinations of $n$, $m$, and $v$. This by itself is enough
to give us $\mathbf m(5)\ge 32$. For $n=6$ and above, we do not get an immediate
improvement. For $n=6$, we need to check the values of $v$ from  $39$ 
to $42$. It turns out that our upper bound for the probability that 
greedy coloring fails is less than $1$ for all cases except $L_1=L_2=9$
and $L_{12}\ge 2$ (we use
$L_A$ to denote the number of edges that contain $A$ as a subset, so,
for example $L_1=l_1+l_{12}+l_{13}+l_{123}$, which equals the degree of vertex
$1$). 
The sum of all vertex degrees is $6\cdot 63=378$ and
every vertex degree is at least $9$, amounting for a total of at least
$39\cdot 9=351$. This only leaves
 room for at most $23$ vertices with a degree greater than 
$9$. So, at least $16$ vertices must have a degree equal to $9$. If we pick one
of those, call it $1$, the sum of the numbers of pairs involving it is 45. 
Every such pair must occur at least once, leaving at most $6$ pairs that 
can occur more than once. Thus,
 we can find another vertex of degree $9$, call it
$2$, such that the pair $(1,2)$ has only one occurrence. Choosing these two 
vertices as the marked ones, we get $L_1=L_2=9$ and $L_{12}=1$. For these
parameters, the greedy algorithm succeeds with positive probability, and
we arrive at $\mathbf m(6)\ge 64$. Similar
arguments work for eliminating $n=8, m=262$ and $n=9,m=537$, so we can 
summarize
\begin{theorem}
We have the lower bounds
\[\mathbf m(5)\ge 32, \mathbf m(6)\ge 64, \mathbf m(7)\ge 128, \mathbf m(8)\ge 263, \mathbf m(9)\ge 538.\]
\end{theorem}

\section{The case $v=2n+1$}
The case $v=2n+1$ is interesting, because it is the smallest number of
vertices for which $\mathbf m(n,v)$ is not known for all values of $n$, and
because of its close relation with other covering questions. In particular,
de Vries\cite{deVries1977} showed that the Goldberg-Russell lower bound
\[v(n,2n+1)\ge\frac{\genfrac(){0pt}1{2n+1}{n}}{n+2}\]
is attained if and only if the associated hypergraph is a $(n-1,n,2n+1)$ 
Steiner system, i.e., if every $n-1$-subset of $V$ is contained in exactly
one of the hyperedges. It is worth noting that this lower bound is
$C_{n+1}/2$, where $C_n$ denotes the $n$-th Catalan number. $C_n$ is known
to be odd iff $n$ is of the form $n=2^k-1$\cite{deutschsagan},
 so we find that for $n=2^k-2$
the lower bound $C_{n+1}/2$ is not an integer, so it cannot be attained. 
De Vries' result, combined with the well-known fact that a $(3,4,v)$-Steiner
system exists iff $v\equiv 4,6\pmod 6$ \cite{hanani}, even gives that
$\mathbf m(n,2n+1)>C_{n+1}/2$ if $n\not\equiv 3,5\pmod 6$. For particular values
of $n$, we can push this a bit further:
\begin{theorem}
For $n=2^k-2$, $k\ge 3$
\[\mathbf m(n,2n+1)\ge (C_{n+1}+3)/2,\]
and for $n=2^r(2^{2k}+1)-2$, $r,k\ge 1$, in particular 
for $n=2^{2k+1}$, $k\ge 1$,
\[\mathbf m(n,2n+1)\ge C_{n+1}/2+3.\]
\end{theorem}

For the proof, we exploit our analysis for the greedy algorithm with one and 
two fixed vertices. We let 
\[m_0=\frac{C_{n+1}}{2}=\frac{(2n+1)!}{n!(n+2)!}\]
 and 
\[L_0=\frac{m_0n}{2n+1}=\frac{(2n)!}{(n-1)!(n+2)!}
=\frac{1}{n-1}\genfrac(){0pt}1{2n}{n-2}=
=\frac{1}{n+2}\genfrac(){0pt}1{2n}{n-1}. \]
The last two terms show that the
 denominator in the representation of $L_0$ as a reduced fraction is a 
common divisor of $n-2$ and $n+1$. By our assumptions, $n-1$ is not a multiple
of $3$, so $L_0$ is an integer. 

Let us now look at a hypergraph with $m=m_0+x$ edges.
A critical vertex can only occur in position $n$, $n+1$, or $n+2$. 

If we lock a vertex with minimal degree $L_1$ (we will call it ``1'' in
the future) in position $n+1$, we can 
bound the probability that $n$ is critical by the probability that it is
the last in some edge, which in turn is bounded by
\[\frac{m-L_1}{\genfrac(){0pt}1{2n}{n}},\]
and the same bound applies to the probability that $n+2$ is critical. 
Similarly, the probability that $n+1$ is critical is bounded by
\[\frac{L_1\genfrac(){0pt}1{n}{n-1}}{\genfrac(){0pt}1{2n}{n-1}}.\]
From these, we obtain the upper bound
\[\frac{ (n!)^2(2m+L_1(n-1))}{(2n)!}\]
for the probability that there is a critical vertex.

The pigeonhole principle implies
\[L_1\le\frac{mn}{2n+1}=\frac{(m_0+x)n}{2n+1}=L_0+\frac{xn}{2n+1}.\]
For $x\le 2$ this upper bound is less then $L_0+1$, so we obtain
$L_1\le L_0$. For $L_1<L_0$, our upper bound for the probability
of a critical vertex is at most
\[1+ \frac{(n!)^2(2x+1-n)}{(2n)!}<1,\]
so we are left with the case $L_1=L_0$. We lock the vertex with the
second largest degree $L_2$ (call it ``2'')
in position $n+2$. We let $L_{12}$ denote
the number of edges containing both vertices 1 and 2.
There are $l_1=L_1-L_{12}$ edges containing vertex 1 but not 2,
$l_2=L_2-L_{12}$ containing vertex 2 alone, and $l_0=m-L_1-L_2+L_{12}$
edges containing neither 1 nor 2. We get upper bounds
\[\frac{l_0}{\genfrac (){0pt}1{2n-1}{n}}\]
and
\[\frac{l_2}{\genfrac (){0pt}1{2n-1}{n-1}}\]
for the probabilities that $n$ resp.\ $n+2$ are critical.
For $n+1$, we have the two upper bounds
\[\frac{l_1n}{\genfrac(){0pt}1{2n-1}{n-1}}\]
and
\[\frac{l_1}{\genfrac(){0pt}1{2n-1}{n-1}}+
\frac{L_{12}(n-1)}{\genfrac(){0pt}1{2n-1}{n-2}}.\]
Thus, the probability that there is a critical vertex is at most
\[\frac{n!(n-1)!(l_0+l_2+\min(nl_1, l_1+L_{12}(n+1)))}{(2n-1)!}=\]\[
\frac{n!(n-1)!(2m+(n-1)L_1-2n|L_{12}-L_{12}^{0}|)}{2(2n-1)!}=\]\[
1+\frac{n!(n-1)!(x-n|L_{12}-L_{12}^{0}|)}{(2n-1)!},\]
where we have put
\[L_{12}^{0}=\frac{L_0(n-1)}{2n}=\frac{m_0(n-1)}{2(2n+1)}=
\frac{(2n-1)!}{(n-2)!(n+2)!}.\]
By theorem 2.1 in \cite{deutschsagan}, the multiplicity of $2$ as a factor
of $C_n$ is one less than the number of digits $1$ in the binary representation
of $n+1$. From this, we can conclude that $L_{12}^0$ is not an integer, but
$4L_{12}^0$ is. Thus $|L_{12}-L_{12}^0|\ge 1/4$, and we get an upper bound
less than $1$ for the probability that a critical vertex exists.

\section{Conclusion}
We were able to obtain improved lower bounds for $\mathbf m(n)$ and $\mathbf m(n,v)$ for
selected values of $n$ and $v$, in particular we could improve the lower
bound $\mathbf m(5)\ge 29$ from \cite{shanni} to $\mathbf m(5)\ge 32$. 

We hope that our method can be extended to obtain an improvement of the
asymptotic lower bound. To this end, we would need to let the number of
fixed vertices go to infinity as $n$ increases. In order to still obtain
sufficiently small upper bounds for the probabilities of criticality,
this needs to be combined with a tight control of the joint occurrences of
the selected vertices.

\section{Acknowledgements}
We are very gratefully to Prof.\ Cherkashin for his valuable remarks and suggestions.

\end{document}